\documentclass{elsart}
\usepackage{amssymb, amsmath, url,graphicx}
\journal{Comptes Rendus Mathematique}

\def\3{\subset }
\def\4{\subseteq }
\def\<{\left<}
\def\>{\right>}

\def\bit{\begin{itemize}}
\def\eit{\end{itemize}}
\def\3{\subset }
\def\4{\subseteq }
\def\ov{\overline}

\def\0{\leqno}

\def\barr{\begin{array}}
\def\earr{\end{array}}

\def\Z{{\rlap{$\kern2pt{\rm Z}$}{\rm Z}\,}}

\begin{document}
\begin{frontmatter}

\title{A characterization of generalized\\ quaternion 2-groups}
\author{Marius T\u arn\u auceanu}
\ead{tarnauc@uaic.ro}
\address{Faculty of Mathematics,
\textquotedblleft Al.I. Cuza\textquotedblright\ University, Ia\c si, Romania}

\begin{abstract}
The goal of this note is to give a characterization of generalized
quaternion 2-groups by using their posets of cyclic subgroups.
\end{abstract}

\end{frontmatter}

\section{Introduction}

Let $G$ be a finite group and $L(G)$ be the subgroup lattice of
$G$. The starting point for our discussion in given by the paper
\cite{1}, where the proper nontrivial subgroups $H$ of $G$ with
the property that

$$\text{for every}\hspace{1mm} X\in L(G)\hspace{1mm} \text{we have}\hspace{1mm} X\leq H\hspace{1mm} \text{or}\hspace{1mm} H\leq X$$

\noindent have been studied. Such a subgroup is called a
\textit{breaking point} for the lattice $L(G)$. Clearly, if $L(G)$
is a chain (i.e. $G$ is a cyclic $p$-group), then all proper
nontrivial subgroups $H$ of $G$ are breaking points. On the other
hand, we remark that the above concept can naturally be extended
to other remarkable posets of subgroups of $G$ (and also to
arbitrary posets). One of them is the poset of cyclic subgroups of
$G$, denoted usually by $C(G)$. The study of the existence and of
the uniqueness of breaking points in $C(G)$ constitutes the
purpose of this paper.

Most of our notation is standard and will usually not be repeated
here. Elementary concepts and results on group theory can be found
in \cite{2} and \cite{4}. For subgroup lattice notions we refer
the reader to \cite{3} and \cite{5}.

We mention that by a \textit{generalized quaternion {\rm 2}-group}
we mean a group of order $2^n$ for some natural number $n\geq 3$,
defined by the presentation
$$Q_{2^n}=\langle a,b \mid a^{2^{n-2}}= b^2, a^{2^{n-1}}=1, b^{-1}ab=a^{-1}\rangle.$$We also
recall that these groups are the unique finite noncyclic
$p$-groups all of whose abelian subgroups are cyclic, or
equivalently the unique finite noncyclic $p$-groups possessing
exactly one subgroup of order $p$ (see (4.4) of \cite{4}, II).
Obviously, this result shows that the subgroup of order 2 of
$Q_{2^n}$, namely $\langle a^{2^{n-2}}\rangle$, is the unique
breaking point of $C(Q_{2^n})$.

Our main theorem proves that generalized quaternion 2-groups
exhaust all finite noncyclic groups whose posets of cyclic
subgroups have breaking points.

\noindent{\bf Theorem 1.1.} {\it Let $G$ be a finite group. Then
$C(G)$ possesses breaking points if and only if $G$ is either a
cyclic $p$-group of order at least $p^2$ or a generalized
quaternion {\rm 2}-group.}

\section{The proof of Theorem 1.1}

We observe first that the above theorem can be easily proved in
the particular case of $p$-groups.

\noindent{\bf Lemma 2.1.} {\it Let $G$ be a finite $p$-group. Then
$C(G)$ possesses breaking points if and only if $G$ is either a
cyclic $p$-group of order at least $p^2$ or a generalized
quaternion {\rm 2}-group.}

\noindent{\bf Proof.} Suppose that $G$ is not cyclic and let $H$
be a breaking point of $C(G)$. Then all minimal subgroups $M_1$,
$M_2$, ..., $M_k$ of $G$ are contained in $H$. If $k\geq 2$, then
we infer that $H$ is not cyclic, a contradiction. So, we have
$k=1$, that is $G$ has a unique subgroup of order $p$. This
implies that $G$ is a generalized quaternion 2-group, according to
the result mentioned in Section 1.

The converse implication is obvious, completing the proof.
\hfill\rule{1,5mm}{1,5mm}

We are now able to give a proof of Theorem 1.1.

\noindent{\bf Proof.} Suppose that the poset $C(G)$ of cyclic
subgroups of a finite group $G$ possesses a breaking point, say
$H$.

In the following we shall focus on proving that $G$ must
necessarily be a $p$-group. By the way of contradiction, assume
that the order of $G$ has at least two distinct prime divisors.
Clearly, the same thing can be also said about the order of $H$.
Let $p\in\pi(G)$ and $K$ be a cyclic $p$-subgroup of $G$. Since
$H$ is not a $p$-subgroup, we infer that $K\subseteq H$. In other
words, $H$ contains any cyclic $p$-subgroup of $G$ and
consequently any $p$-element of $G$. This implies that all Sylow
$p$-subgroups of $G$ are contained in $H$. Then $H=G$, a
contradiction.

Hence $G$ is a $p$-group, for some prime $p$, and now the
conclusion follows from Lemma 2.1.
\hfill\rule{1,5mm}{1,5mm}

By the above results we also infer that, given a finite group $G$,
the poset $C(G)$ possesses a \textit{unique} breaking point if and
only if $G$ is either a cyclic $p$-group of order $p^2$ or a
generalized quaternion 2-group. In other words, the following
corollary holds.

\noindent{\bf Corollary 2.2.} {\it The generalized quaternion {\rm
2}-groups are the unique finite noncyclic groups whose posets of
cyclic subgroups have exactly one breaking point.}

Finally, we indicate a natural generalization of our study,
suggested by the reviewers of the paper. Let $G$ be a finite group
and denote by $$\ov{C}(G)=\{\hspace{0,5mm}[H]\mid H\in C(G)\}$$the
set of conjugacy classes of cyclic subgroups of $G$. Mention that
$\ov{C}(G)$ is also a poset under the ordering relation
$$[H_1]\leq [H_2]\mbox{ if and only if } H_1\subseteq H_2^g, \mbox{ for some } g\in
G.$$Take a breaking point $[H]$ of $\ov{C}(G)$. Then $H\in C(G)$
satisfies the following condition: for any cyclic subgroup $C$ of
$G$, some conjugate of $C$ in $G$ contains or is contained in $H$.
Clearly, this is weaker than the condition that $H$ be a breaking
point of $C(G)$. We remark that for a finite $p$-group $G$ it is
sufficient to guarantee the uniqueness of a subgroup of order $p$
in $G$. In other words, Lemma 2.1 also holds if we replace $C(G)$
with $\ov{C}(G)$. In the general case, that is for
\textit{arbitrary} finite groups $G$, the problem of
characterizing the existence and the uniqueness of breaking points
of $\ov{C}(G)$ remains still open.

\section{Conclusions and further research}

All previous results show that the concept of breaking point in
some posets of subgroups of a (finite) group $G$ can constitute an
important aspect of subgroup lattice theory. Clearly, its study
started in \cite{1} for lattices of subgroups and continued in the
present paper for posets of cyclic subgroups can successfully be
extended to other significant lattices/posets associated to $G$
(as the lattice of normal/subnormal/characteristic/solitary
subgroups of $G$ or the poset of centralizers/conjugacy classes of
elements (subgroups) of $G$). Studying the breaking points of
arbitrary posets (not necessarily connected with a group $G$)
seems to be also very interesting. These will surely constitute
the subject of some further research.
\bigskip

{\bf Acknowledgements.} The author is grateful to the reviewers
for their remarks which improve the previous version of the paper.

\end{document}